\documentclass[11pt]{article}
\usepackage[utf8]{inputenc}
\usepackage[T1]{fontenc}
\DeclareUnicodeCharacter{211D}{\ensuremath{\mathbb{R}}}
\DeclareUnicodeCharacter{211A}{\ensuremath{\mathbb{Q}}}
\DeclareUnicodeCharacter{2102}{\ensuremath{\mathbb{C}}}
\DeclareUnicodeCharacter{2124}{\ensuremath{\mathbb{Z}}}
\DeclareUnicodeCharacter{00D7}{\ensuremath{\times}}
\DeclareUnicodeCharacter{2200}{\ensuremath{\forall}}
\DeclareUnicodeCharacter{2203}{\ensuremath{\exists}}
\DeclareUnicodeCharacter{2208}{\ensuremath{\in}}
\DeclareUnicodeCharacter{00AC}{\ensuremath{\neg}}
\DeclareUnicodeCharacter{2192}{\ensuremath{\to}}
\DeclareUnicodeCharacter{2022}{\ensuremath{\bullet}}
\DeclareUnicodeCharacter{1D65}{\ensuremath{_v}}
\DeclareUnicodeCharacter{2B1D}{\ensuremath{\cdot}}
\DeclareUnicodeCharacter{1D40}{\ensuremath{^{\mathsf T}}}
\DeclareUnicodeCharacter{27E8}{\ensuremath{\langle}}
\DeclareUnicodeCharacter{27E9}{\ensuremath{\rangle}}
\DeclareUnicodeCharacter{03C6}{\ensuremath{\varphi}}
\DeclareUnicodeCharacter{221A}{\ensuremath{\surd}}
\DeclareUnicodeCharacter{2264}{\ensuremath{\le}}
\DeclareUnicodeCharacter{2227}{\ensuremath{\wedge}}
\DeclareUnicodeCharacter{00B3}{\ensuremath{^3}}
\usepackage{amsmath,amssymb,amsthm,mathtools}
\usepackage{upquote}
\usepackage{booktabs}
\usepackage[margin=1in]{geometry}
\usepackage[colorlinks=true,linkcolor=blue,citecolor=blue,urlcolor=blue]{hyperref}

\newtheorem{theorem}{Theorem}[section]
\newtheorem{lemma}[theorem]{Lemma}
\newtheorem{proposition}[theorem]{Proposition}

\theoremstyle{definition}

\theoremstyle{remark}

\newcommand{\R}{\mathbb{R}}
\newcommand{\Q}{\mathbb{Q}}
\newcommand{\C}{\mathbb{C}}
\newcommand{\Z}{\mathbb{Z}}
\newcommand{\Fr}{\operatorname{Fr}}
\newcommand{\diag}{\operatorname{diag}}
\newcommand{\re}{\operatorname{Re}}
\newcommand{\im}{\operatorname{Im}}
\newcommand{\T}{\mathsf T}
\DeclarePairedDelimiter{\abs}{\lvert}{\rvert}
\DeclarePairedDelimiter{\norm}{\lVert}{\rVert}

\title{No regular icosahedron has a rational point on every face:\\ a proof checked in Lean}
\author{Mikhail Levit}
\date{November 8, 2026}

\begin{document}

\maketitle

\begin{abstract}
We prove that no regular icosahedron in $\R^3$ has a point with rational coordinates on each of its
twenty closed faces, in every position, size and orientation, reflections included. This answers a
question of I.~Bogdanov on MathOverflow. A field automorphism of $\C$ that sends $\sqrt5$ to
$-\sqrt5$ reduces the problem to the following statement: for every complex orthogonal $3\times3$
matrix $T$ some real linear combination of $\re T$ and $\im T$ separates ten intervals determined by
the faces. We prove this statement by a computer-assisted case analysis over rational
parametrizations of rotation frames, in which every inequality is an exact polynomial bound over
$\Z[\varphi]$, $\varphi=(1+\sqrt5)/2$; no floating-point arithmetic and no Lipschitz estimates are
used. The whole argument, including the $48$ certificates, is formalized in the Lean~4 proof
assistant with Mathlib and checked by the Lean kernel; the main theorem depends only on the
standard axioms. The mathematical argument and the formalization were produced by AI systems
(OpenAI Codex and Anthropic Claude) under the direction of the author; see the statement at the end of
the paper.
\end{abstract}

\section{Introduction}

Can a regular icosahedron contain a rational point on each of its faces? The question was asked
by Ilya Bogdanov on MathOverflow in August 2023 \cite{Bogdanov} and has remained without an answer
there. For the other Platonic solids the answer is positive and easy: regular tetrahedra, cubes and
octahedra may have all vertices rational, and eight vertices of a regular dodecahedron form a cube
whose edges are diagonals of the twelve pentagonal faces \cite{Bogdanov}. A regular icosahedron
cannot have all vertices rational, since there is no regular pentagon in an integer lattice of any
dimension; but this alone says nothing about one rational point per face.

\begin{theorem}\label{thm:main}
Let $I\subset\R^3$ be a regular icosahedron. Then at least one of the twenty closed triangular faces
of $I$ contains no point of $\Q^3$.
\end{theorem}

Points on edges and vertices are allowed; a face is the closed triangle. Closedness matters: in the
comments to \cite{Bogdanov} J.~Rouse gives a regular icosahedron and six rational points that
together lie on all twenty face \emph{planes}, so that the analogue of the question for planes has a
positive answer. The planar analogue of the theorem is known: there is no regular pentagon with a
rational point on each edge \cite{Pentagon}. That proof is algebraic and works even for rational
points on the lines through the edges; M.~Bolan extended it to all regular $n$-gons, which have a
rational point on each edge only for $n=3,4,8$ \cite{Pentagon}. By Rouse's example, no argument
of this kind can work for the icosahedron: a proof has to use the boundaries of the faces.

Our proof has two parts. The first part (Section~\ref{sec:red}) is algebraic. The planes of the faces are
defined over $K=\Q(\sqrt5)$; applying an automorphism $\tau$ of $\C$ with $\tau(\sqrt5)=-\sqrt5$ to the
equations of the faces and the rational points turns a hypothetical solution into a complex
orthogonal matrix $T\in SO(3,\C)$ such that, for every real $\alpha,\beta$, ten real intervals built
from the matrix $\alpha\re T+\beta\im T$ and the ten pairs of opposite faces have a common point.
The second part (Sections~\ref{sec:sep} and~\ref{sec:cert}) shows that this never happens: for every
$T\in SO(3,\C)$ some combination separates the intervals. Linear algebra reduces this to a statement
about a compact parameter space, which is then checked on finitely many cells by exact arithmetic.

The whole proof is formalized in Lean~4 \cite{Lean4} with Mathlib \cite{Mathlib} and checked by the
Lean kernel (Section~\ref{sec:formal}). The computational part is carried out by the kernel itself
through \texttt{decide +kernel} on $48$ certificates; \texttt{native\_decide} and other ways to trust
compiled code are not used. The only parts that must be trusted are the Lean kernel and the definitions
in the statement: $58$ lines that define the canonical icosahedron, its faces and rational points.

\section{Reduction to complex orthogonal matrices}\label{sec:red}

Let $\varphi=(1+\sqrt5)/2$, $\varphi'=1-\varphi=(1-\sqrt5)/2$, $K=\Q(\sqrt5)$, and let $x\mapsto x'$
denote the nontrivial automorphism of $K$. The canonical icosahedron $I_0$ has the twelve
vertices obtained from $(0,\pm1,\pm\varphi)$ by cyclic permutations of the coordinates, and edge
length~$2$. Every regular icosahedron is the image of $I_0$ under $x\mapsto t+sRx$ with $t\in\R^3$,
$s>0$ and $R$ orthogonal; this is how Theorem~\ref{thm:main} is stated in Lean.

The twenty faces $F_0,\dots,F_{19}$ of $I_0$ come in ten pairs of opposite faces; we fix
representatives $F_0,\dots,F_9$, so that the face opposite to $F_i$ is $-F_i$. Each face $F_j$ lies in a
plane $n_j\cdot x=\varphi^2$ with $n_j\in K^3$: the eight normals are $(\pm1,\pm1,\pm1)$ and the
twelve others are cyclic permutations of $(0,\pm\varphi,\pm(\varphi-1))$. We write $n_j'$ for the
conjugate normal and $v_{i0},v_{i1},v_{i2}$ for the vertices of $F_i$.

\begin{lemma}[centering]\label{lem:center}
Suppose that both $t+sRF_i$ and the opposite face $t+sR(-F_i)$ contain rational points. Then there is
$m\in F_i$ with $sRm\in\Q^3$.
\end{lemma}

\begin{proof}
If $t+sRx=p$ and $t-sRy=q$ with $x,y\in F_i$ and $p,q\in\Q^3$, then $m=(x+y)/2\in F_i$ by convexity
and $sRm=(p-q)/2$.
\end{proof}

\begin{lemma}\label{lem:aut}
There is a field automorphism $\tau$ of $\C$ with $\tau(\sqrt5)=-\sqrt5$, that is,
$\tau(\varphi)=\varphi'$.
\end{lemma}

\begin{proof}
Extend the conjugation of $K$ to a purely transcendental extension of $K$ in $\C$ generated by a
transcendence basis, and then to its algebraic closure $\C$.
\end{proof}

For a real matrix $A$ and $0\le i\le9$, $0\le k\le2$ put
\[
  f_{ik}(A)=n_i'^{\T}A\,v_{ik},\qquad J_i(A)=\bigl[\min_k f_{ik}(A),\ \max_k f_{ik}(A)\bigr].
\]
Since $v\mapsto n_i'^{\T}Av$ is linear, $J_i(A)$ is the set of its values on the triangle $F_i$. We say
that $A$ \emph{separates} if $J_0(A),\dots,J_9(A)$ have no common point. Clearly, if $A$
separates, then so does $cA$ for every $c>0$.

\begin{proposition}\label{prop:red}
If every face of $t+sRI_0$ contains a rational point, then the matrix $T=(R^\tau)^{\T}R$, where
$\tau$ is applied entrywise, satisfies $T^{\T}T=I$ and $\det T=1$, and no matrix
$\alpha\re T+\beta\im T$ with $\alpha,\beta\in\R$ separates.
\end{proposition}

\begin{proof}
By Lemma~\ref{lem:center} there are $m_i\in F_i$ with $p_i=sRm_i\in\Q^3$. Since $R^{\T}R=I$ and
$n_i\cdot m_i=\varphi^2$, we have $(Rn_i)\cdot p_i=s\varphi^2$. Apply $\tau$: it fixes $p_i$, maps the
entries of $n_i$ to those of $n_i'$ and $\varphi$ to $\varphi'$, so $(R^\tau n_i')\cdot p_i=\tau(s)\varphi'^2$.
Substituting $p_i=sRm_i$ gives
\[
  n_i'^{\T}\,T\,m_i=\kappa\,\varphi'^2,\qquad \kappa=\tau(s)/s\in\C,\qquad i=0,\dots,9.
\]
Since $n_i'$, $m_i$ and $\varphi'^2$ are real, for $A=\alpha\re T+\beta\im T$ we get
$n_i'^{\T}Am_i=(\alpha\re\kappa+\beta\im\kappa)\varphi'^2$, the same number for all $i$, and it lies in
every $J_i(A)$. Finally, $T^{\T}T=R^{\T}\tau(RR^{\T})R=I$ and $\det T=\tau(\det R)\det R=(\det R)^2=1$,
because $\det R=\pm1$ is fixed by $\tau$.
\end{proof}

Hence Theorem~\ref{thm:main} follows from the next theorem, which is the main work.

\begin{theorem}[separation]\label{thm:sep}
For every $T\in M_3(\C)$ with $T^{\T}T=I$ and $\det T=1$ there are $\alpha,\beta\in\R$ such that
$\alpha\re T+\beta\im T$ separates.
\end{theorem}

\section{Proof of the separation theorem}\label{sec:sep}

\subsection{Chart frames}
For $u,v\in\R$ let $Q(u,v)$ be the rotation matrix of the quaternion $1+u\,\mathbf i+v\,\mathbf j$
without normalization:
\[
  Q(u,v)=\begin{pmatrix}1+u^2-v^2&2uv&2v\\2uv&1-u^2+v^2&-2u\\-2v&2u&1-u^2-v^2\end{pmatrix},\qquad
  \mu(u,v)=1+u^2+v^2 .
\]
Then $Q^{\T}Q=QQ^{\T}=\mu^2I$ and $\det Q=\mu^3$, so $Q/\mu\in SO(3)$. Let $P_0,\dots,P_5$ be the
signed permutation matrices with columns $e_i,\sigma e_j,\sigma e_k$, where $(i,j,k)$ is a cyclic shift
of $(1,2,3)$ and $\sigma=\pm1$; they satisfy $P_c^{\T}P_c=I$, $\det P_c=1$. The \emph{chart frame} is
$\Fr_c(u,v)=P_cQ(u,v)$. Its third column is $P_c(2v,-2u,1-u^2-v^2)^{\T}$, which is $\mu$ times
a unit vector.

\begin{lemma}[six charts]\label{lem:cover}
For every unit vector $d\in\R^3$ there are $c\in\{0,\dots,5\}$ and $u,v$ with $\abs u,\abs v\le1/2$
and $\Fr_c(u,v)e_3=\mu(u,v)\,d$.
\end{lemma}

\begin{proof}
Choose $c$ so that $P_c^{\T}d=(a_1,a_2,a_3)$ has $a_3\ge\abs{a_1},\abs{a_2}$ (the coordinate of $d$ of
largest absolute value, with its sign). Put $u=-a_2/(1+a_3)$, $v=a_1/(1+a_3)$. From
$a_1^2+a_2^2=(1-a_3)(1+a_3)$ one gets $\mu=2/(1+a_3)$ and $(2v,-2u,1-u^2-v^2)=\mu(a_1,a_2,a_3)$;
moreover $\abs u,\abs v\le a_3/(1+a_3)\le1/2$.
\end{proof}

Let $J=\left(\begin{smallmatrix}0&-1\\1&0\end{smallmatrix}\right)$ and, for real numbers,
\[
  M_r(\sigma,t)=\diag\bigl(\sigma((1-t^2)I_2+2tJ),\ 1+t^2\bigr),\qquad
  M_c(x,y,z)=\diag\bigl(xI_2+yJ,\ z\bigr).
\]
The computational part of the proof is the following theorem; its certificates are described in
Section~\ref{sec:cert}.

\begin{theorem}[certificates]\label{thm:cert}
\begin{enumerate}
\item[\textup{(R)}] For every chart $c$, every sign $\sigma=\pm1$ and all $\abs u,\abs v\le1/2$,
  $\abs t\le1$ the matrix $\Fr_c(u,v)M_r(\sigma,t)$ separates.
\item[\textup{(C)}] For every pair of charts $c_1,c_2$ and all $u,v,u',v'\in[-1/2,1/2]$ there is a triple
  $(X,Y,D)$ from Table~\ref{tab:phases}, depending on $(u,v,u',v')$, such that for all
  $z\in[-1,1]$ the matrix $\Fr_{c_1}(u,v)\,M_c(X,Y,Dz)\,\Fr_{c_2}(u',v')^{\T}$ separates.
\end{enumerate}
\end{theorem}

\begin{table}[ht]
\centering
\begin{tabular}{c|ccccccccc}
\toprule
$X$ & $1$ & $99$ & $12$ & $91$ & $119$ & $161$ & $5$ & $20$ & $0$\\
$Y$ & $0$ & $20$ & $5$ & $60$ & $120$ & $240$ & $12$ & $99$ & $1$\\
$D$ & $1$ & $101$ & $13$ & $109$ & $169$ & $289$ & $13$ & $101$ & $1$\\
\bottomrule
\end{tabular}
\caption{The sixteen phases $(\pm X,Y,D)$, $X^2+Y^2=D^2$: the nine listed triples and the seven
triples $(-X,Y,D)$ with $X\neq0,1$. They are directions in the upper half plane.}
\label{tab:phases}
\end{table}

\subsection{The real case}

Write $X=\re T$, $Y=\im T$. From $T^{\T}T=I$:
\begin{equation}\label{eq:parts}
  X^{\T}X=I+Y^{\T}Y,\qquad X^{\T}Y+Y^{\T}X=0 .
\end{equation}

\begin{lemma}[frame decomposition]\label{lem:frame}
Let $T^{\T}T=I$, $\det T=1$, let $b\in\R^3$ satisfy $Yb=0$, and let the chart frames
$F_1=\Fr_{c_1}(u,v)$, $F_2=\Fr_{c_2}(u',v')$ satisfy $F_1e_3=\mu\,Xb$ and $F_2e_3=\mu'b$, where
$\mu=\mu(u,v)$, $\mu'=\mu(u',v')$. Then there are $p,r\in\C$ with $p^2+r^2=1$ and
\[
  F_1^{\T}\,T\,F_2=\mu\mu'\,\diag(pI_2+rJ,\ 1).
\]
\end{lemma}

\begin{proof}
Let $M=F_1^{\T}TF_2$ and $k=\mu\mu'$. Then $M^{\T}M=F_2^{\T}T^{\T}(F_1F_1^{\T})TF_2=k^2I$ and
$\det M=k^3$. Since $Yb=0$, $TF_2e_3=\mu'Tb=\mu'Xb=(\mu'/\mu)F_1e_3$, hence
$Me_3=(\mu'/\mu)F_1^{\T}F_1e_3=ke_3$. Then $M^{\T}e_3=k^{-1}M^{\T}Me_3=ke_3$, so the third row and
column of $M$ are $(0,0,k)$. The upper left block $P$ satisfies $P^{\T}P=k^2I_2$ and
$k\det P=k^3$; for $P=\left(\begin{smallmatrix}a&b'\\c&d\end{smallmatrix}\right)$ this gives $d=a$,
$b'=-c$, and $p=a/k$, $r=c/k$.
\end{proof}

\begin{lemma}\label{lem:real}
If $\im T=0$, then $\re T=\gamma\,\Fr_c(u,v)M_r(\sigma,t)$ for some $\gamma>0$, chart $c$, sign $\sigma$
and $\abs u,\abs v\le1/2$, $\abs t\le1$.
\end{lemma}

\begin{proof}
By \eqref{eq:parts} $X$ is orthogonal, so $Xe_3$ is a unit vector. Take $c,u,v$ from
Lemma~\ref{lem:cover} for $d=Xe_3$ and apply Lemma~\ref{lem:frame} with $b=e_3$ and
$F_2=\Fr_0(0,0)=I$. Since $T=X$ is real, so are $p,r$, and $p^2+r^2=1$. If $p\ge0$ put $\sigma=1$,
$t=r/(1+p)$; otherwise $\sigma=-1$, $t=r/(p-1)$. In both cases $\abs t\le1$ and
$M_r(\sigma,t)=(1+t^2)\diag(pI_2+rJ,1)$ (the half-angle substitution). Since
$F_1F_1^{\T}=\mu^2I$, we get $X=\mu^{-2}F_1(F_1^{\T}X)=(\mu(1+t^2))^{-1}F_1M_r(\sigma,t)$.
\end{proof}

\subsection{The non-real case}

\begin{lemma}\label{lem:ker}
There is a unit vector $b\in\R^3$ with $Yb=0$; for it $\norm{Xb}=1$.
\end{lemma}

\begin{proof}
By \eqref{eq:parts}, $\norm{Xw}^2=\norm w^2+\norm{Yw}^2$ for every $w$, so $X$ is invertible. From
$X^{\T}Y=-Y^{\T}X$ we get $\det X\det Y=-\det Y\det X$, hence $\det Y=0$. Take $b$ in the kernel of
$Y$; then $\norm{Xb}^2=\norm b^2=1$.
\end{proof}

\begin{proposition}\label{prop:span}
Let $\im T\neq0$. There are charts $c_1,c_2$ and $u,v,u',v'\in[-1/2,1/2]$ such that for every
$(X_0,Y_0,D)$ with $D>0$ and $X_0^2+Y_0^2=D^2$ there are $\alpha,\beta\in\R$, $\gamma>0$ and $z\in[-1,1]$
with
\[
  \alpha\re T+\beta\im T=\gamma\,\Fr_{c_1}(u,v)\,M_c(X_0,Y_0,Dz)\,\Fr_{c_2}(u',v')^{\T}.
\]
\end{proposition}

\begin{proof}
Take $b$ from Lemma~\ref{lem:ker}, frames $F_1$ for $Xb$ and $F_2$ for $b$ from
Lemma~\ref{lem:cover}, and $p=p_1+ip_2$, $r=r_1+ir_2$ from Lemma~\ref{lem:frame}. Taking real
and imaginary parts,
\[
  F_1^{\T}XF_2=\mu\mu'M_c(p_1,r_1,1),\qquad F_1^{\T}YF_2=\mu\mu'M_c(p_2,r_2,0),
\]
and $p^2+r^2=1$ gives $a=p_1^2+r_1^2=1+b''$, $b''=p_2^2+r_2^2$, and $p_1p_2+r_1r_2=0$. If $b''=0$,
then $F_1^{\T}YF_2=0$ and $Y=(\mu\mu')^{-2}F_1(F_1^{\T}YF_2)F_2^{\T}=0$; so $b''>0$ and $a\ge1$. Put
\[
  \alpha=(X_0p_1+Y_0r_1)\,b'',\qquad \beta=(X_0p_2+Y_0r_2)\,a,\qquad \lambda=ab'',\qquad
  z=\alpha/(\lambda D).
\]
Using $p_1p_2+r_1r_2=0$ one checks $\alpha p_1+\beta p_2=\lambda X_0$ and $\alpha r_1+\beta r_2=\lambda Y_0$,
so $F_1^{\T}(\alpha X+\beta Y)F_2=\mu\mu'\lambda\,M_c(X_0,Y_0,Dz)$. By Cauchy--Schwarz
$\abs{X_0p_1+Y_0r_1}\le D\sqrt a\le Da$, hence $\abs z\le1$. Multiplying by $F_1$ on the left and by
$F_2^{\T}$ on the right and using $F_iF_i^{\T}=\mu_i^2I$ gives the claim with
$\gamma=\lambda/(\mu\mu')$.
\end{proof}

\begin{proof}[Proof of Theorem~\ref{thm:sep}]
If $\im T=0$, take $\alpha=1$, $\beta=0$: by Lemma~\ref{lem:real} and Theorem~\ref{thm:cert}(R),
$\re T$ is a positive multiple of a separating matrix. If $\im T\neq0$, take the charts and parameters
of Proposition~\ref{prop:span}; Theorem~\ref{thm:cert}(C) gives a phase $(X_0,Y_0,D)$ that separates
for all $z\in[-1,1]$, and Proposition~\ref{prop:span} with this phase gives $\alpha,\beta$.
\end{proof}

The frames in Proposition~\ref{prop:span} do not depend on the phase, while the phase chosen by
Theorem~\ref{thm:cert}(C) depends on the frames; this order of quantifiers is what the proof needs.

\section{The certificates}\label{sec:cert}

\subsection{Exact arithmetic}
All checks are computations with natural numbers that the Lean kernel evaluates directly; its
natural numbers are arbitrary precision and addition, multiplication and comparison are
implemented by GMP. An element $a+b\varphi$ of $\Z[\varphi]$ is stored as four naturals
$(a^+,a^-,b^+,b^-)$. Its sign is decided exactly: with $w=2a+b$, $a+b\varphi\ge0$ iff
$w+b\sqrt5\ge0$, which is decided by comparing $w^2$ with $5b^2$ when $w$ and $b$ have
opposite signs. The checker uses only \texttt{Nat.add}, \texttt{Nat.mul}, \texttt{Nat.sub},
\texttt{Nat.ble}, structures, \texttt{cond} and explicit recursion \texttt{Nat.rec}, so that
\texttt{decide +kernel} can reduce it.

\subsection{Cells and polynomial bounds}
A cell of a chart square at level $m$ and index $j$ is $u=(2j+1-2^m+s)/2^{m+1}$, $s\in[-1,1]$; at level
$0$ it is the whole $[-1/2,1/2]$. A $t$-interval at level $m$ is $t=(2j+1-2^m+s)/2^m$. With
$N=2^{m+1}$ the entries of $N^2Q(u,v)$ are quadratic polynomials in the two cell coordinates with
integer coefficients, and so are the entries of $2^{2m}M_r(\sigma,t)$ in the coordinate of $t$. Hence,
for $A=\Fr_{c_1}M\Fr_{c_2}^{\T}$,
\[
  f_{ik}(A)=\bigl(\Fr_{c_1}^{\T}n_i'\bigr)\cdot\bigl(M\Fr_{c_2}^{\T}v_{ik}\bigr)
  =\text{(positive constant)}\cdot a_i\cdot w_{ik},
\]
where $a_i$ and $w_{ik}$ are vectors of polynomials of degree at most $2$ in two variables over
$\Z[\varphi]$; similarly in the real case with $w_{ik}=M_rv_{ik}$, a polynomial in one variable. The
positive constant does not affect separation. Write $a=a^0+\tilde a$ and $w=w^0+\tilde w$ (constant
terms and the rest). On $[-1,1]^2$ the polynomial $c_1s+c_2t+c_3s^2+c_4st+c_5t^2$ is at least
$\min(c_3,0)+\min(c_5,0)-\abs{c_1}-\abs{c_2}-\abs{c_4}$, and
$\abs{\tilde a\cdot\tilde w}\le\sum_{l}\norm{\tilde a_l}_1\norm{\tilde w_l}_1$, where
$\norm{\cdot}_1$ is the sum of the absolute values of the coefficients. This gives an exact lower
bound of $f_{i_1k}-f_{i_2k'}$ on the whole cell. A pair $(i_1,i_2)$ \emph{separates the cell} if these
nine lower bounds, $k,k'\in\{0,1,2\}$, are positive; then $J_{i_1}(A)$ lies to the right of $J_{i_2}(A)$
for every $A$ of the cell, and $A$ separates. The soundness of this bound is proved in Lean.

\subsection{The real case: search in the kernel}
For (R) no certificate data is needed. The function \texttt{rsearch} takes a cell, chooses at its
centre $i_1$ maximizing $\min_kf_{ik}$ and $i_2$ minimizing $\max_kf_{ik}$, and checks whether
$(i_1,i_2)$ separates the cell; if not, it splits the $(u,v)$ square into four or the $t$-interval
into two (the coarser one) and recurses, with a recursion bound of $40$. The kernel runs this search
for the six charts and two signs, twelve runs in total, and the results are the twelve theorems
\texttt{rsearch c $\sigma$ 40 $\langle0,0,0,0,0\rangle$ = true}.

\subsection{The complex case: certificate trees}
For (C) the search over $(u,v,u',v')$ is guided by a tree, stored as one natural number whose
hexadecimal digits are read from the least significant end: $0$ means a split of the $(u,v)$ or the
$(u',v')$ square (the coarser one) into four; $1\,n\,i_1\,i_2$ is a leaf where phase $n$ and pair
$(i_1,i_2)$ separate the cell at $z=-1$ and at $z=1$; $2\,n\,i_1\,i_2\,j_1\,j_2$ is a leaf where
$(i_1,i_2)$ separates at $z=-1$ and $z=0$ and $(j_1,j_2)$ at $z=0$ and $z=1$. The forms $f_{ik}$ are
affine in $z$, so separation by a fixed pair at both ends of a $z$-interval implies it on the whole
interval. The trees are produced by a generator (\texttt{Gen.lean}) that runs the same functions
compiled and checks its output with them; the generator is not part of the proof. Table~\ref{tab:stats}
gives the sizes.

\begin{table}[ht]
\centering
\begin{tabular}{lr}
\toprule
pairs of charts & $36$\\
leaves with one pair & $14\,472$\\
leaves with two pairs & $4\,896$\\
splits & $6\,444$\\
maximal total level $m_1+m_2$ & $6$\\
hexadecimal digits in all trees & $93\,708$\\
\bottomrule
\end{tabular}
\caption{The certificate trees of the complex case.}
\label{tab:stats}
\end{table}

\subsection{Resources}
On an Intel Core i5-12500 with $32$~GB of memory the kernel checks one real certificate in
$0.5$--$2$ minutes using up to $7$~GB, and one complex certificate in $1.5$--$3$ minutes using $8$--$13$~GB.
Checking the complex certificates one at a time, all $48$ take about $1.5$ hours. Building them with the
default parallelism of Lake exhausts $32$~GB of memory, so the repository contains scripts that build
them sequentially.

\section{The formalization}\label{sec:formal}

The Lean project \cite{Repo} uses Lean~4 v4.34.0-rc2 and Mathlib v4.34.0-rc2. Table~\ref{tab:files}
lists the files. The main theorem reads
\begin{verbatim}
theorem no_rational_point_on_every_face (t : Fin 3 → ℝ) (s : ℝ) (hs : 0 < s)
    (R : Matrix (Fin 3) (Fin 3) ℝ) (hR : R ∈ Matrix.orthogonalGroup (Fin 3) ℝ) :
    ¬ ∀ j : Fin 20, HasRatPoint t s R j
\end{verbatim}
and its meaning depends only on the definitions in \texttt{Defs.lean}:
\begin{verbatim}
noncomputable def φ : ℝ := (1 + Real.sqrt 5) / 2
noncomputable def ev (x : ℤ × ℤ) : ℝ := x.1 + x.2 * φ
noncomputable def vertex (i : Fin 12) : Fin 3 → ℝ := fun c => ev (vtxCode i c)
noncomputable def face (j : Fin 20) : Set (Fin 3 → ℝ) :=
  convexHull ℝ (Set.range fun k => vertex (faceIdx j k))
def IsRat (x : Fin 3 → ℝ) : Prop := ∃ q : Fin 3 → ℚ, ∀ c, x c = q c
def HasRatPoint (t : Fin 3 → ℝ) (s : ℝ) (R : Matrix (Fin 3) (Fin 3) ℝ)
    (j : Fin 20) : Prop :=
  ∃ x ∈ face j, IsRat (t + s • (R *ᵥ x))
\end{verbatim}
where \texttt{vtxCode} lists the twelve vertices $(0,\pm1,\pm\varphi)$, $(\pm1,\pm\varphi,0)$,
$(\pm\varphi,0,\pm1)$ as pairs $(a,b)$ for $a+b\varphi$, and \texttt{faceIdx} lists the twenty faces as
triples of vertex indices. The orthogonal group includes reflections. A short script in the
repository checks that \texttt{faceIdx} consists exactly of the twenty triples of pairwise adjacent
vertices, that is, of the faces. The command \texttt{\#print axioms} gives
\begin{verbatim}
'Icosahedron.no_rational_point_on_every_face' depends on axioms:
  [propext, Classical.choice, Quot.sound]
\end{verbatim}

\begin{table}[ht]
\centering
\begin{tabular}{lrl}
\toprule
file & lines & content\\
\midrule
\texttt{Defs} & $58$ & the statement: icosahedron, faces, rational points\\
\texttt{Geometry} & $132$ & face planes, opposite faces, conjugate normals\\
\texttt{Embedding} & $88$ & Lemma~\ref{lem:aut}\\
\texttt{Reduction} & $165$ & Lemma~\ref{lem:center}, Proposition~\ref{prop:red}\\
\texttt{KCheck} & $487$ & the checker: $\Z[\varphi]$, polynomials, searches, phases\\
\texttt{KSound} & $349$ & soundness of the arithmetic and of the cell bound\\
\texttt{KGeom} & $535$ & chart frames; soundness of both searches\\
\texttt{Frames} & $646$ & Lemmas~\ref{lem:cover}--\ref{lem:ker}, Proposition~\ref{prop:span}\\
\texttt{CertR/*}, \texttt{CertC/*} & & $48$ certificates (\texttt{decide +kernel}), $110$~KB\\
\texttt{CertAll} & $131$ & collects the certificates\\
\texttt{Separation} & $38$ & Theorem~\ref{thm:sep}\\
\texttt{Main} & $21$ & Theorem~\ref{thm:main}\\
\bottomrule
\end{tabular}
\caption{Files of the Lean library \texttt{Icosahedron}.}
\label{tab:files}
\end{table}

\section*{Use of AI tools}

The author posed the problem, directed the work and checked the results. The mathematical argument
(the reduction of Section~\ref{sec:red} and the decomposition of Section~\ref{sec:sep}) was developed
by OpenAI Codex, and the Lean formalization, including the
kernel-checkable certificates, their generator and the final form of the decomposition, was developed by
Anthropic Claude, both under the direction of the author; the text of this paper was drafted with
Anthropic Claude. The author takes full responsibility for the content. The correctness of
Theorem~\ref{thm:main} does not depend on trust in these systems: the proof is checked by the Lean
kernel, and what must be trusted is the kernel and the $58$ lines of definitions quoted in
Section~\ref{sec:formal}.

\end{document}